\documentclass{article}
\usepackage{latexsym,amssymb,amsmath}
\begin{document}
\begin{center}\Large{On a non-combinatorial definition of Stirling numbers}\end{center}
 \begin{center}Milan Janji\'c\end{center}
 \begin{center}Faculty of Natural Sciences and Mathematics, \end{center}
\begin{center} Banja Luka, Republic of Srpska,
Bosnia and Herzegovina.\end{center}
\begin{center}e-mail: agnus@blic.net\end{center}
\begin{center}\textbf{Introduction}\end{center}

In Combinatorics Stirling  numbers may be defined in several way.
One such definition is given in [1], where an extensive
consideration of Stirling numbers is presented.

In this paper an alternative definition of Stirling numbers of
both kind is given. Namely, Stirling numbers of the first kind appear in the closed formula for the $n$-th derivative of $f(\ln x)$.
 In the same way Stirling numbers of the second kind appear in the formula for the $n$-th derivative of $f(e^x).$
Here $f(x)$ is  an  arbitrary smooth real function.

 This facts allow us to define
 Stirling numbers within the frame of differential calculus. These definitions may be interesting because arbitrary functions appear in them.
Choosing suitable function we may obtain different properties of Stirling numbers by the use of derivatives only.

Using simple properties of derivatives we obtain  here three important properties of Stirling numbers.
First  are so called two terms recurrence relations, from which one can easily derive the combinatorial meaning of Stirling numbers.

Next we obtain expansion formulas  of powers into falling factorials, and vise versa.
These expansions usually serve as the definitions of Stirling
numbers, as in [1].

Finally, we obtain the exponential generating functions for Stirling and Bell numbers.

As a by product the closed formulas for the $n$-th derivative of the
functions $f(e^x)$ and $f(\ln x)$ are obtained.

\begin{center}\textbf{Stirling numbers of the first kind}\end{center}

Suppose that $f(x)\in C^{\infty}(0,+\infty).$
We shall calculate $[f(\ln x)]^{(n)}$ by $x.$ Denoting $\ln x=t$ we
have

$$ [f(\ln x)]'=\frac{f'(t)}{x}$$$$[f(\ln x)]''=-\frac{f'(t)}{x^2}+\frac{f''(t)}{x^2},$$
$$[f(\ln x)]'''=2\frac{f'(t)}{x^3}-3\frac{f''(t)}{x^3}+\frac{f'''(t)}{x^3}.$$
Since
\begin{equation}\label{r2}\bigg[\frac{f^{(k)}(t)}{x^m}\bigg]'_x=\frac{-mf^{(k)}(t)+f^{(k+1)}(t)}{x^{m+1}}\end{equation} we
may obtain $[f(\ln x)]^{(n)}$ in the form

\begin{equation}\label{sb1}f(\ln x)]^{(n)}=\frac{1}{x^n}\sum_{k=1}^ns(n,k)f^{(k)}(t).\end{equation}

It is clear that $s(n,k)$ are  integers not depending on $f(x).$

\noindent\textbf{Definition 1.} \textit{The integers
$s(n,k),\;(n\geq 1,\;k=1,\ldots,n)$ are called  Stirling numbers
of the first  kind. We also define $s(0,0)=1$ and $s(n,k)=0,$ for
$n<k$.}

The first result that we shall derive is two terms recurrence relation [1, Th. 8.7 (a)].

It holds

$$s(n,1)=(-1)^{n-1} (n-1)!,\;s(n,n)=1,\;$$$$s(n+1,k)=s(n,k-1)-ns(n,k),\;(k=2,3,\ldots,n).$$

To prove this we take derivative in
$$x^nf^{(n)}(\ln x)=\sum_{k=1}^ns(n,k)f^{(k)}(t)$$
and conclude that

$$nx^{n-1}f^{(n)}(\ln x)+x^{n-1}f^{(n+1)}(\ln x)=\frac 1x \sum_{k=1}^ns(n,k)f^{(k+1)}(t),$$ that is,
$$x^{n+1}f^{(n+1)}(\ln x)=\sum_{k=1}^ns(n,k)f^{(k+1)}(t)-n\sum_{k=1}^ns(n,k)f^{(k)}(t)=$$$$=
-ns(n,1)f'(t)+\sum_{k=2}^n[s(n,k-1)-ns(n,k)]f^{(k)}(t)+s(n,n)f^{(n+1)}(t).$$

We thus obtain
$$s(n+1,1)=-ns(n,1),\;s(n+1,n+1)=s(n,n),$$$$
s(n+1,k)=s(n,k-1)-ns(n,k),\;(k=2,3,\ldots,n).$$

From this the desired equations  follows easily.

We next prove the formula of expansion of falling factorials into powers, that usually serves as the definition, as in  [1, Def. 8.1].

 For each
$n\geq 1$ and each $\alpha\in\mathbb{R}$ holds
\begin{equation}\label{vo}\alpha(\alpha-1)\cdots (\alpha-n+1)=\sum_{k=1}^ns(n,k)\alpha^k.\end{equation}

To prove this it is enough to take $f(x)=x^{\alpha}$ in (\ref{sb1}).

Replacing $x$ by $1+x$ in (\ref{sb1}) we get
\begin{equation}\label{ee}(1+x)x^nf^{(n)}[\ln (1+x)]=\sum_{k=1}^ns(n,k)f^{(k)}(t),\end{equation}
where $t=\ln(1+x).$

Take a fixed $k$ and consider the function $f(t)=\frac{t^k}{k!}.$
We have $f[\ln(1+x)]=\frac{[\ln(1+x)]^k}{k!}.$ It is clear that
$f^{(k)}(1)=k!,$ and $f^{(n)}(1)=0,\;(n\not=k).$ Replacing this in
(\ref{ee}) we obtain the following [1, Th. 8.3 (a)]:

For
$n\geq k$ holds
$$\left\{\frac{[\ln(1+x)]^k}{k!}\right\}^{(n)}_{{x=0}}=s(n,k).$$
In other wards, $$\frac{[\ln(1+x)]^k}{k!}$$ is the exponential
generation function for Stirling numbers of the first kind.

\begin{center}\textbf{Stirling numbers of the second kind and Bell numbers}\end{center}

 Calculating
 derivatives of
$f(e^x)$ by $x$ we obtain
$$ [f(e^x)]'=tf'(t)$$$$[f(e^x)]''=tf'(t)+t^2f''(t),$$
$$[f(e^x)]'''=tf'(t)+3t^2f''(t)+t^3f'''(t),$$
where $t=e^x.$

According to the fact that
\begin{equation}\label{r1}[t^kf^{(k)}(t)]'_x=kt^kf^{(k)}(t)+t^{k+1}f^{(k+1)}(t)\end{equation} we
may obtain $[f(e^x)]^{(n)}$ in the form

\begin{equation}\label{sb2}[f(e^x)]^{(n)}=\sum_{k=1}^nS(n,k)t^kf^{(k)}(t).\end{equation}

It is clear that $S(n,k)$ are  integers not depending on $f(x).$

\noindent\textbf{Definition 2.} \textit{The integers
$S(n,k),\;(n\geq 1,\;k=1,\ldots,n)$ are called  Stirling numbers
of the second  kind. We also define $S(0,0)=1$ and $S(n,k)=0,$ for
$n<k$.}

Taking derivative in (\ref{sb2}) and applying (\ref{r1}) we obtain
$$[f(e^x)]^{(n+1)}=\sum_{k=1}^nS(n,k)\big[kt^kf^{(k)}(t)+t^{k+1}f^{(k+1)}(t)\big]=$$$$=
\sum_{k=1}^nS(n,k)kt^kf^{(k)}(t)+\sum_{k=2}^{n+1}S(n,k-1)kt^kf^{(k)}(t)=$$$$=
S(n,1)f'(t)+\sum_{k=1}^n[kS(n,k)+S(n,k-1)]t^kf^{(k)}(t)+S(n,n)t^{n+1}f^{(n+1)}(t).$$
From this we conclude that
$$S(n+1,1)=S(n,1),\;S(n+1,n+1)=S(n,n),$$$$S(n+1,k)=kS(n,k)+S(n,k-1),\;(k=2,3,...,n).$$
We have thus proved the following [1, Th. 8.7 (b)]:

Stirling numbers of the second kind fulfill the following
recurrence relations
$$S(n,1)=S(n,n)=1,$$$$S(n,k)=S(n-1,k-1)+kS(n-1,k),\;(k=2,3,\ldots,n).$$

If we take specially  $f(x)=x^\alpha,\;(\alpha\in \mathbb R)$ we
have $f(e^x)=e^{\alpha x}$ which implies
$[f(e^x)]^{(n)}=\alpha^ne^{\alpha x}.$ On the other hand, we have
$f^{(k)}(t)=\alpha(\alpha-1)\cdots(\alpha-k+1)x^{\alpha-k}.$ By
replacing this into (\ref{sb2}) we obtained the following [1, Def. 8.1]:

For each
$\alpha\in\mathbb R$ holds
$$\alpha^n=\sum_{k=1}^nS(n,k)\alpha(\alpha-1)\cdots(\alpha-k+1).$$

Consider the function $f(x)=\frac{(x-1)^k}{k!}.$ It is clear that
$$f^{(k)}(x)=1,\;f^{(n)}(x)=0,\;(n>k).$$
In the case $n<k$ holds $f^{(n)}(1)=0.$ It follows that
$$\left[\sum_{m=1}^nS(n,m)t^mf^{(m)}(t)\right]_{t=1}=S(n,k).$$

Using (\ref{sb2}) we obtain the following [1, Th. 8.3 (b)]:

For
$n\geq k$ holds
$$S(n,k)=\bigg[\frac{(e^x-1)^k}{k!}\bigg]^{(n)}_{x=0}.$$
Equivalently, $$\frac{(e^x-1)^k}{k!}$$ is the exponential generating
function for Stirling numbers of the second kind.

Taking $f(x)=e^x$ in the equation (\ref{sb2}) gives the following
assertions for Bell polynomials $B_n(x)$ and Bell numbers $B_n.$

 For $n=1,2,\ldots$ hold
$$e^x\big(e^{e^x}\big)^{(n)}=e^{e^x}B_n(t),\;(t=e^x)$$
$$e^{-1}\big(e^{e^x}\big)^{(n)}_{x=0}=B_n.$$
The last equation is equivalent to the fact that  $e^{-1}e^{e^x}$
is the exponential generating function for Bell numbers.

\vspace{1cm}

\begin{center}\textbf{Reference}\end{center}

\noindent [1] C. A. Charalambides. \textit{Enumerative
combinatorics}, Chapman\&Hall/CRC, 2002.

\end{document}